\documentclass[12pt]{article}
\usepackage[pdftex]{color, graphicx}
\usepackage{setspace}
\usepackage[utf8]{inputenc}
\usepackage{amsmath}
\usepackage[pdftex,bookmarks,colorlinks]{hyperref}
\usepackage{latexsym}
\usepackage{amsfonts}
\usepackage{amssymb}
\usepackage{amsthm}
\usepackage{paralist}
\usepackage{mathrsfs}
\usepackage[super]{nth}

\begin{document}\newtheorem{thm}{Theorem}
\newtheorem{cor}[thm]{Corollary}
\newtheorem{lem}{Lemma}
\theoremstyle{remark}\newtheorem{rem}{Remark}
\theoremstyle{definition}\newtheorem{defn}{Definition}
\title{Variation and oscillation inequalities for operator averages on a complex Hilbert space}
\author{Sakin Demir\\
Agri Ibrahim Cecen University\\ 
Faculty of Education\\
\nth{3} Floor, Office C-42\\
04100 A\u{g}r{\i}, Turkey\\
e-mail: sakin.demir@gmail.com
}
\date{06.21.2025}
\maketitle
\renewcommand{\thefootnote}{}



\renewcommand{\thefootnote}{\arabic{footnote}}
\setcounter{footnote}{0}

\begin{abstract} Let $\mathcal{H}$ be a complex Hilbert space and $T:\mathcal{H}\to \mathcal{H}$  be a contraction. Let
$$A_nf=\frac{1}{n}\sum_{j=1}^nT^jf$$
for $f\in \mathcal{H}$.  Let $(n_k)$  be a lacunary sequence, then  there exists a constant $C_1>0$ such that
$$\sum_{k=1}^\infty\|A_{n_{k+1}}f-A_{n_k}f\|_{\mathcal{H}}\leq C_1\|f\|_{\mathcal{H}}$$
for all $f\in \mathcal{H}$.\\
\indent
Let $(n_k)$  be a lacunary  sequence, and let $\mathbb{N}$ be the set of natural numbers. Then there exists a constant $C_2>0$ such that 
$$\sum_{k=1}^\infty\sup_{\substack{n_k\leq m< n_{k+1}\\m\in \mathbb{N}}}\|A_m(T)f-A_{n_k}(T)f\|_{\mathcal{H}}\leq C_2\|f\|_{\mathcal{H}}$$
for all $f\in \mathcal{H}$.
\end{abstract}
\section{Introduction}
\indent

Let $\mathcal{H}$ be a complex Hilbert space, and $T:\mathcal{H}\to \mathcal{H}$. The adjoint $T^\ast$ of $T$ is the unique mapping $T^\ast:\mathcal{H}\to \mathcal{H}$ such that 
$$\langle Tf , g \rangle =\langle f , T^\ast g\rangle $$
for all $f,g\in \mathcal{H}$. $T$ is called self-adjoint if $T=T^\ast$. $T$ is called an isometry if $T^\ast T=I$ (indentity on $\mathcal{H}$), and a coisometry if its adjoint is an isometry, i.e., $TT^\ast=I$. It is a unitary operator if it is an isometry and a coisometry; that is, if it is an invertible transformation with $T^{-1}=T^\ast$ or equivalently, an invertible isometry which means a surjective isometry. \\
An operator $T:\mathcal{H}\to \mathcal{H}$ is called a contraction  if $\|Tf\|_{\mathcal{H}}\leq \|f\|_{\mathcal{H}}$ for every  $f\in \mathcal{H}$.\\
We say that $f,g\in \mathcal{H}$ are orthogonal if $\langle f, g\rangle =0$. In this case we write $f\perp g$. Two subsets $A$ and $B$ of $\mathcal{H}$ are orthogonal  if every element of $A$ is orthogonal to every element of $B$. In this case we write $A\perp B$. The orthogonal complement of a set $A$ is the set $A^{\perp}$ defined as
$$A^{\perp}=\left\{f\in \mathcal{H}: \langle f,  g\rangle =0 \; \text{for every } g\in A\right\},$$
A function $F:X\to X$ of a set $X$ into itself is idempotent if $F^2=F$, where $F^2$ stands for the composition of $F$ with itself. A projection $E$ is an idempotent linear transformation $E:\mathcal{X}\to\mathcal{X}$ of a linear space $\mathcal{X}$ into itself. If $E$ is a projection, then so $I-E$, which is referred to as the complementary projection of $E$, and their null spaces and ranges are related as follows:
$$\mathcal{R}(E)=\mathcal{N}(I-E)\;\;\text{and}\;\; \mathcal{N}(E)=\mathcal{R}(I-E).$$
An orthogonal projection $E:\mathcal{X}\to\mathcal{X}$ on an inner product space $\mathcal{X}$ is a projection such that $\mathcal{R}(E)\perp\mathcal{N}(E)$.\\
\indent
Let $\mathcal{H}$ be a complex Hilbert space and $T:\mathcal{H}\to \mathcal{H}$ be a contraction. Define the average
$$A_nf=\frac{1}{n}\sum_{j=1}^nT^jf$$
for any  $f\in \mathcal{H}$.\\
\indent
The following result has been proved in R. L. Jones~{\it{et al}}~\cite{liftweb}:
\begin{thm}\label{jlrsblt}Let $(n_k)$ be an increasing sequence of positive integers, and let $T$ be a contraction on a complex Hilbert space $\mathcal{H}$. Then
$$\left(\sum_{k=1}^\infty\|A_{n_{k+1}}(T)f-A_{n_k}(T)f\|_{\mathcal{H}}^2\right)^{1/2}\leq 25 \|f\|_{\mathcal{H}}$$
for all $f\in \mathcal{H}$.
\end{thm}
Four years later M.~Lifshits, M.~Weber~\cite{liftweb} showed in the following result that the constant of  Theorem~\ref{jlrsblt} can be improved:
\begin{thm}\label{jrblt}Let $(n_k)$ be an increasing sequence of positive integers, and let $T$ be a contraction on a complex Hilbert space $\mathcal{H}$. Then
$$\sum_{k=1}^\infty\|A_{n_{k+1}}(T)f-A_{n_k}(T)f\|_{\mathcal{H}}^2\leq 6\pi \|f\|^2_{\mathcal{H}}$$
for all $f\in \mathcal{H}$.
\end{thm}
\indent
The following reverse versions of this result have been obtained by the author in S. Demir~\cite{sdemrinq}:
\begin{thm}\label{sdemrinq}Let $U:\mathcal{H}\to \mathcal{H}$ be a unitary operator on a Hilbert space $\mathcal{H}$ and $A_nf=\frac{1}{n}\sum_{j=1}^nU^jf$ for all $f\in \mathcal{H}$. Suppose that $(n_k)$ is a lacunary sequence with no non-trivial common divisor, then there exists a positive constant $C$ such that
$$\|f\|_{\mathcal{H}}\leq C\left(\sum_{k=1}^\infty \|A_{n_{k+1}}f-A_{n_k}f\|_{\mathcal{H}}^2\right)^{1/2}$$
for all $f\in  \mathcal{H}$ with $\int f=0$.
\end{thm}
\begin{thm}\label{sdcont}
Let $T$ be a contraction on a Hilbert space $\mathcal{H}$ and let $A_n(T)f=\frac{1}{n}\sum_{j=1}^nT^jf$ for all $f\in \mathcal{H}$. Suppose that  $(n_k)$ is a lacunary sequence with no non-trivial common divisor, then there exist a Hilbert space $K$ containing $\mathcal{H}$ as a closed subspace, and an orthogonal projection $P:K\to  \mathcal{H}$ such that
$$\|P\|_{\mathcal{H}}\|f\|_{\mathcal{H}}\leq C\left(\sum_{k=1}^\infty \|A_{n_{k+1}}(T)f-A_{n_k}(T)f\|_{\mathcal{H}}^2\right)^{1/2}.$$
for all $f\in \mathcal{H}$ with $\int f=0$, where $C$ is a positive constant.
\end{thm}
The following result has been obtained by the author in a recent work in S. Demir~\cite{sdemiraoi}:
\begin{thm}\label{aoi}
 Let $T$ be a contraction on a complex Hilbert space $\mathcal{H}$,  and for $f\in \mathcal{H}$ define
$$A_n(T)f=\frac{1}{n}\sum_{j=1}^nT^jf.$$ 
Let $(n_k)$ be an increasing sequence and $M$ be any sequence. We prove that there exists a positive constant $C$ such that
$$\left(\sum_{k=1}^\infty\sup_{\substack{n_k\leq m< n_{k+1}\\m\in M}}\|A_m(T)f-A_{n_k}(T)f\|_{\mathcal{H}}^2\right)^{1/2}\leq C\|f\|_{\mathcal{H}}$$
for all $f\in \mathcal{H}$.
\end{thm}
\indent
Oscillation and variation inequalities have become an indispensable tool in harmonic analysis, operator theory, ergodic theory and martingale theory to study problems of pointwise convergence. They also have become a subject of interest in their of right.\\
\indent
It is clear that
\begin{equation}
\left(\sum_{k=1}^\infty\|A_{n_{k+1}}(T)f-A_{n_k}(T)f\|_{\mathcal{H}}^2\right)^{1/2}\leq \sum_{k=1}^\infty\|A_{n_{k+1}}(T)f-A_{n_k}(T)f\|_{\mathcal{H}}
\end{equation}
and
\begin{equation}
\left(\sum_{k=1}^\infty\sup_{\substack{n_k\leq m< n_{k+1}\\m\in M}}\|A_m(T)f-A_{n_k}(T)f\|_{\mathcal{H}}^2\right)^{1/2}\;\;\;\;\;\;\;\;\;\;\;\;\;\;\;\;\;\;\;\;\;\;\;\;\;\;\;\;\;\;\;\;\;\;\;\;\;\;\;\;
\end{equation}
\begin{equation*}
\;\;\;\;\;\;\;\;\;\;\;\;\;\;\;\;\;\;\;\;\;\;\;\;\;\;\;\;\;\;\;\;\;\;\;\;\;\;\;\;\;\;\;\;\;\;\;\;\;\;\; \leq \sum_{k=1}^\infty\sup_{\substack{n_k\leq m< n_{k+1}\\m\in M}}\|A_m(T)f-A_{n_k}(T)f\|_{\mathcal{H}}
\end{equation*}
with the above settings.\\
\indent
It is the goal l of this research to obtain optimal conditions to be able to control the right hand sides of (1) and (2), since none of the above mentioned results is sufficient to deal with this issue in general.\\
\indent
Note that the Hilbert spaces are all complex Hilbert spaces in this research when it is not stated otherwise.
\section{Results}
\indent

Let $\mathcal{H}$ be a complex Hilbert space and $T:\mathcal{H}\to \mathcal{H}$ be a contraction. Define the average
$$A_nf=\frac{1}{n}\sum_{j=1}^nT^jf$$
for any  $f\in \mathcal{H}$.\\
It is known by the spectral theorem for a contraction, for instance see \cite{liftweb}, that there exists a measure $\mu$ on $[-\pi ,\pi )$ such that $\mu [-\pi , \pi )=\|f\|_{\mathcal{H}}^2$ and
$$\|A_m(f)-A_n(f)\|_{\mathcal{H}}\leq \|a_m(\theta )-a_n(\theta )\|_{2,\mu}$$
where
$$a_n(\theta )=\frac{1}{n}\left(1+e^{i\theta}+e^{i2\theta}+\dots +e^{i(n-1)\theta}\right)=\frac{e^{in\theta}-1}{n(e^{i\theta}-1)}$$
with equality when $T$ is an isometry.\\
\indent
Recall that a sequence $(n_k)$ is called lacunary if there exists a real number $\alpha >1$ such that
$$\frac{n_{k+1}}{n_k}\geq \alpha$$
for all $k=1,2,3,\dots $\\
\indent
The following variation inequality is our first result:
\begin{thm}\label{varineq} Let $T$ be a contraction on a complex Hilbert space $\mathcal{H}$, and for $f\in \mathcal{H}$ define
$$A_nf=\frac{1}{n}\sum_{j=1}^nT^jf.$$
If  $(n_k)$ is a lacunary sequence , then there exists a constant $C>0$ such that 
$$\sum_{k=1}^\infty\|A_{n_{k+1}}f-A_{n_k}f\|_{\mathcal{H}}\leq C\|f\|_{\mathcal{H}}$$
for all $f\in \mathcal{H}$.
\end{thm}
\begin{proof} It suffices to show that
$$\sum_{k=1}^N\|A_{n_{k+1}}f-A_{n_k}f\|_{\mathcal{H}}\leq C\|f\|_{\mathcal{H}}$$
for all $N\geq 1$. \\
\indent
Because of the spectral theorem theorem for a contraction on a complex Hilbert space (see U.~Krengel~\cite{ukgl}, p. 94) there exists a probability measure $\mu $ on $[-\pi , \pi ]$ such that
$$\langle T^nf, T^mf\rangle =\int_{-\pi}^{\pi}e^{i(m-n)\theta}\, d\mu (\theta )$$
for all $n,m\in \mathbb{Z}$, where $\langle \cdot , \cdot \rangle$ denotes the inner product on $\mathcal{H}$. Therefore when we set
$$a_n(\theta )=\frac{e^{in\theta}-1}{n(e^{i\theta}-1)}$$
it follows that  whenever $n,m\geq 1$ we have
\begin{align*}
\|A_m(f)-A_n(f)\|_{\mathcal{H}}^2&\leq \|a_m(\theta )-a_n(\theta )\|_{2,\mu}.
\end{align*}
We now have
\begin{align*}
\|A_m(f)-A_n(f)\|_{\mathcal{H}}&\leq \left(\int_{-\pi}^{\pi}|a_m(\theta )-a_n(\theta )|^2\, d\mu (\theta )\right)^{1/2}.
\end{align*}
It is thus clear that
\begin{align*}
\sum _{k=1}^N\|A_{n_{k+1}}(f)-A_{n_k}(f)\|_{\mathcal{H}}&\leq \sum _{k=1}^N\left(\int_{-\pi}^{\pi}|a_{n_{k+1}}(\theta )-a_{n_k}(\theta )|^2\, d\mu (\theta )\right)^{1/2}.
\end{align*}
We can assume without the loss of generality that $\theta \in (0, \pi ]$.  Since $|e^{i\theta}-1|\geq\frac{1}{4}\theta$ for all $\theta\in (0, \pi  ]$, we have
\begin{align*}
|a_{n_{k+1}}(\theta )-a_{n_k}(\theta )|&=\left|\frac{1}{n_{k+1}}\left(\frac{e^{in_{k+1}\theta}-1}{e^{i\theta}-1}\right)-\frac{1}{n_k}\left(\frac{e^{in_k\theta}-1}{e^{i\theta}-1}\right)\right|\\
&\leq 4 \left|\frac{e^{in_{k+1}\theta}-1}{n_{k+1}\theta}-\frac{e^{in_k\theta}-1}{n_k\theta}\right|
\end{align*}
and therefore it suffices to control the summation
$$I(\theta )=\sum _{k=1}^N\left(\int_{-\pi}^{\pi}\left|\frac{e^{in_{k+1}\theta}-1}{n_{k+1}\theta}-\frac{e^{in_k\theta}-1}{n_k\theta}\right|^2\, d\mu (\theta )\right)^{1/2}.$$
First fix $\theta$ and let $k_0$ be the first $k$ such that $\theta n_k\geq 1$, and let
$$I_1(\theta )=\sum_{k:\theta n_k<1}\left(\int_{-\pi}^{\pi}\left|\frac{e^{in_{k+1}\theta}-1}{n_{k+1}\theta}-\frac{e^{in_k\theta}-1}{n_k\theta}\right|^2\, d\mu (\theta )\right)^{1/2},$$
$$I_2(\theta )=\sum_{k:\theta n_k\geq 1}\left(\int_{-\pi}^{\pi}\left|\frac{e^{in_{k+1}\theta}-1}{n_{k+1}\theta}-\frac{e^{in_k\theta}-1}{n_k\theta}\right|^2\, d\mu (\theta )\right)^{1/2}.$$
Since $I(\theta )=I_1(\theta )+I_2(\theta )$,  it suffices to control $I_1(\theta )$ and $I_2(\theta )$ to prove our theorem.\\
\indent
To control $I_1(\theta )$  first define a function $F:\mathbb{R}\to\mathbb{C}$ as
$$F(r)=\frac{e^{ir\theta}-1}{r\theta}.$$
Then we have
$$F(n)=\frac{e^{in\theta}-1}{n\theta}.$$
By the Mean Value Theorem there exits a constant $\xi\in (n_k, n_{k+1})$ such that
$$|F(n_{k+1})-F(n_k)|= |F^\prime (\xi )| |n_{k+1}-n_k|.$$
We have
\begin{align*}
F^\prime (r)&=\frac{i\theta e^{ir\theta} r\theta-\theta (e^{ir\theta}-1)}{(r\theta )^2}\\
&=\frac{e^{ir\theta }(ir\theta -1)+1}{r^2\theta}\\
&=-\frac{(1-ir\theta )e^{ir\theta }-1}{r^2\theta}.
\end{align*}
Thus we obtain
$$|F^\prime (r)|=\frac{|(1-ir\theta )e^{ir\theta }-1|}{|r^2\theta |}$$
and on the other hand there exists a positive constant $C$ such that
$$|(1-z)e^z-1|\leq C|z|^2$$
for $|z|\leq 1$. Thus for $|r\theta |\leq 1$ then we have $|(1-ir\theta )e^{ir\theta }-1|\leq C|r^2\theta^2|$, and it follows that
$$|F^\prime (r)|\leq C|\theta |.$$
Since  $\xi\theta \leq n_{k+1}\theta <1$, we obtain
\begin{align*}
|F(n_{k+1})-F(n_k)|&= |F^\prime (\xi )| |n_{k+1}-n_k|\\
&\leq C\theta (n_{k+1}-n_k)\\
&\leq C\theta n_{k+1}, \;\;\;\; ({\text{since}}\; \;0 <n_k<n_{k+1}) \\
&= C\theta n_{k+1}\cdot \frac{n_{k_0-1}}{n_{k_0-1}}\\
&=C\frac{n_{k+1}}{n_{k_0-1}}\cdot n_{k_0-1}\theta \\
&<C\frac{n_{k+1}}{n_{k_0-1}}, \;\;\;\; ({\text{since}}\; \;n_{k_0-1}\theta <1).
\end{align*}
Since
$$\frac{n_{k+1}}{n_{k_0-1}}=\frac{n_{k+1}}{n_{k+2}}\cdot\frac{n_{k+2}}{n_{k+3}}\cdot \frac{n_{k+3}}{n_{k+4}}\dots \frac{n_{k_0-2}}{n_{k_0-1}},$$
and becasue of lacunarity we have
$$\frac{n_k}{n_{k+1}}\leq \frac{1}{\alpha }<1$$
for all $k=1,2,3,\dots $ we obtain
$$\frac{n_{k+1}}{n_{k_0-1}}\leq \frac{1}{\alpha^{k_0-2-(k+1)}}=\frac{1}{\alpha^{k_0-k-3}}$$
and obviously we have
$$\sum_{k:\theta n_k<1}\frac{1}{\alpha^{k_0-k-3}}\leq \sum_{k=1}^\infty \frac{1}{\alpha^k}=\frac{1}{1-\frac{1}{\alpha }}=C(\alpha ).$$
Thus we obtain
\begin{align*}
I_1(\theta )&=\sum_{k:n_k\theta <1}\left(\int_{-\pi}^{\pi}|F(n_{k+1}\theta )-F(n_k\theta )|^2\, d\mu (\theta )\right)^{1/2}\\
&<C \sum_{k:n_k\theta <1}\left(\int_{-\pi}^{\pi}\left( C(\alpha )\right)^2\, d\mu (\theta )\right)^{1/2}\\
&\leq C\cdot C(\alpha )\left(\int_{-\pi}^{\pi}\, d\mu (\theta )\right)^{1/2}\\
&= C\cdot C(\alpha )\|f\|_{\mathcal{H}}
\end{align*}
as we wanted. \\
\indent
Let us now control 
$$I_2(\theta )=\sum_{k:\theta n_k\geq 1}\left(\int_{-\pi}^{\pi}\left|\frac{e^{in_{k+1}\theta}-1}{n_{k+1}\theta}-\frac{e^{in_k\theta}-1}{n_k\theta}\right|^2\, d\mu (\theta )\right)^{1/2}.$$
First note that
$$\left|\frac{e^{in_{k+1}\theta}-1}{n_{k+1}\theta}-\frac{e^{in_k\theta}-1}{n_k\theta}\right|\leq \frac{2}{n_{k+1}\theta}+\frac{2}{n_k\theta}\leq \frac{4}{n_k\theta}$$
since $ n_{k+1}\geq n_k$.\\
Let $k_0$ be the first index such that $\theta n_{k_0}\geq 1$. Since $(n_k)$ is lacunary there exists a real number $\alpha >1$ such that 
$$\frac{n_{k+1}}{n_k}\geq \alpha$$
 for all $k=1,2,3,\dots$ Since
$$\frac{n_{k_0}}{n_k}=\frac{n_{k_0}}{n_{k_0+1}}\cdot\frac{n_{k_0+1}}{n_{k_0+2}}\cdot \frac{n_{k_0+2}}{n_{k_0+3}}\dots \frac{n_{k-1}}{n_k},$$
we obtain
$$\frac{n_{k_0}}{n_k}\leq \frac{1}{\alpha^{k-k_0}}.$$
Since $\frac{1}{\theta n_{k_0}}\leq 1$, we obtain
\begin{align*}
\frac{1}{n_k\theta}&=\frac{1}{\theta n_{k_0}}\cdot\frac{n_{k_0}}{n_k}\\
&\leq  \frac{n_{k_0}}{n_k}\\
&\leq  \frac{1}{\alpha^{k-k_0}}.
\end{align*}
It follows that 
\begin{align*}
I_2(\theta )&=\sum_{k:\theta n_k\geq 1}\left(\int_{-\pi}^{\pi}\left|\frac{e^{in_{k+1}\theta}-1}{n_{k+1}\theta}-\frac{e^{in_k\theta}-1}{n_k\theta}\right|^2\, d\mu (\theta )\right)^{1/2}\\
&\leq \sum_{k:\theta n_k\geq 1}\left(\int_{-\pi}^{\pi}\left( \frac{1}{\alpha^{k-k_0}}\right)^2\, d\mu (\theta )\right)^{1/2}\\
&\leq \sum_{k:\theta n_k\geq 1} \frac{1}{\alpha^{k-k_0}}\|f\|_{\mathcal{H}}\\
& \leq C(\alpha )\|f\|_{\mathcal{H}},
\end{align*}
and this finishes our proof.
\end{proof}
\begin{rem}\label{uicont}It is known and easy to prove that a unitary operator $U:\mathcal{H}\to\mathcal{H}$ on a Hilbert space $\mathcal{H}$  is a contraction.
\end{rem}
\begin{proof}Let  $U$ be a unitary operator, and $U^\ast$ be its adjoint operator on $\mathcal{H}$. Then we have $U^\ast U f=If=f$. Thus we have
\begin{align*}
\|Uf\|^2&=\langle Uf , Uf\rangle \\
&=\langle U^\ast Uf , f\rangle \\
&=\langle If , f\rangle \\
&=\langle f , f\rangle .
\end{align*}
Thus we have $\|Uf\|=\sqrt{\langle f, f\rangle}=\|f\|$, and hence we obtain
\begin{align*}
\|U\|&=\sup_{f\neq 0}\frac{\|Uf\|}{\|f\|}\\
&= \sup_{f\neq 0}\frac{\|f\|}{\|f\|}\\
&=1,
\end{align*}
where $\|U\|$ denotes the operator norm of $U$. Thus $U$ is a contraction.
\end{proof}
\begin{cor}\label{ercorvar}Let  $(X,\mathscr{B} ,\mu ,\tau )$ be an ergodic, measure preserving dynamical system and define the ergodic averages
$$\mathcal{A}_nf(x)=\frac{1}{n}\sum_{j=1}^nf\circ\tau^j(x).$$
If  $(n_k)$ is a lacunary sequence, then there exists a positive constant $C$ such that
$$\sum_{k=1}^\infty\|\mathcal{A}_{n_{k+1}}f-\mathcal{A}_{n_k}f\|_2\leq C\|f\|_2$$
for all $f\in L^2(X)$.
\end{cor}
\begin{proof}Note  that $Tf(x)=f\circ\tau (x)$ is a unitary operator,  and thus it is a contraction. Since $L^2(X)$ is a Hilbert space, the result follows from Theorem~\ref{varineq}.
\end{proof}
Our second result is the following oscillation inequality:
\begin{thm}\label{osineq} Let $T$ be a contraction on a Hilbert space $\mathcal{H}$, and for $f\in \mathcal{H}$ define
$$A_nf=\frac{1}{n}\sum_{j=1}^nT^jf.$$
 Let $(n_k)$  be a lacunary sequence, and let $\mathbb{N}$ be the set of natural numbers. Then there exists a constant $C>0$ such that
$$\sum_{k=1}^\infty\sup_{\substack{n_k\leq m< n_{k+1}\\m\in \mathbb{N}}}\|A_mf-A_{n_k}f\|_{\mathcal{H}}\leq C\|f\|_{\mathcal{H}}$$
for all  $f\in \mathcal{H}$.
\end{thm}
\begin{proof} It suffices to show that
$$\sum_{k=1}^N\sup_{\substack{n_k\leq m< n_{k+1}\\m\in \mathbb{N}}}\|A_mf-A_{n_k}f\|_{\mathcal{H}}\leq C\|f\|_{\mathcal{H}}$$
for all $N\geq 1$.\\ 
\indent
We know from the proof of Theorem~\ref{varineq} that
$$\|A_m(f)-A_{n_k}(f)\|_{\mathcal{H}}\leq \left(\int_{-\pi}^{\pi}|a_m(\theta )-a_{n_k}(\theta )|^2\, d\mu (\theta )\right)^{1/2},$$
where 
$$a_n(\theta )=\frac{e^{in\theta}-1}{n(e^{i\theta}-1)}.$$
But since
$$|a_m(\theta )-a_{n_k}(\theta )|^2\, d\mu (\theta )\leq  \sup_{\substack{n_k\leq m< n_{k+1}\\m\in \mathbb{N}}}|a_m(\theta )-a_{n_k}(\theta )|^2\, d\mu (\theta ),$$
we obtain
$$ \int_{-\pi}^{\pi}|a_m(\theta )-a_{n_k}(\theta )|^2\, d\mu (\theta )\leq  \int_{-\pi}^{\pi}\sup_{\substack{n_k\leq m< n_{k+1}\\m\in \mathbb{N}}}|a_m(\theta )-a_{n_k}(\theta )|^2\, d\mu (\theta ).$$
Hence we obtain
$$\sup_{\substack{n_k\leq m< n_{k+1}\\m\in \mathbb{N}}}\|A_m(f)-A_{n_k}(f)\|_{\mathcal{H}}\leq \left(\int_{-\pi}^{\pi}\sup_{\substack{n_k\leq m< n_{k+1}\\m\in \mathbb{N}}}|a_m(\theta )-a_{n_k}(\theta )|^2\, d\mu (\theta )\right)^{1/2}.$$
We then have
\begin{align*}
\sum_{k=1}^N\sup_{\substack{n_k\leq m< n_{k+1}\\m\in \mathbb{N}}}\|A_m(f)-A_{n_k}(f)\|_{\mathcal{H}}\qquad\qquad\qquad\qquad\qquad\qquad\qquad\\
\qquad\qquad\qquad\leq \sum_{k=1}^N \left(\int_{-\pi}^{\pi}\sup_{\substack{n_k\leq m< n_{k+1}\\m\in \mathbb{N}}}|a_m(\theta )-a_{n_k}(\theta )|^2\, d\mu (\theta )\right)^{1/2}.
\end{align*}
Thus in order to prove our theorem it suffices to control the summation
$$J(\theta )=\sum_{k=1}^N \left(\int_{-\pi}^{\pi}\sup_{\substack{n_k\leq m< n_{k+1}\\m\in \mathbb{N}}}|a_m(\theta )-a_{n_k}(\theta )|^2\, d\mu (\theta )\right)^{1/2}.$$
\indent
Let now
$$J_1(\theta )=\sum_{k:n_k\theta <1} \left(\int_{-\pi}^{\pi}\sup_{\substack{n_k\leq m< n_{k+1}\\m\in \mathbb{N}}}|a_m(\theta )-a_{n_k}(\theta )|^2\, d\mu (\theta )\right)^{1/2}$$
and
$$J_2(\theta )=\sum_{k:n_k\theta \geq 1}\left(\int_{-\pi}^{\pi}\sup_{\substack{n_k\leq m< n_{k+1}\\m\in \mathbb{N}}}|a_m(\theta )-a_{n_k}(\theta )|^2\, d\mu (\theta )\right)^{1/2}.$$
We can assume without the loss of generality that $\theta \in (0, \pi ]$. Since $|1-e^{i\theta}|\geq\frac{1}{4}\theta$ for all $\theta\in (0, \pi  ]$ we have
\begin{align*}
|a_m(\theta )-a_{n_k}(\theta )|&=\left|\frac{1}{m}\left(\frac{1-e^{im\theta}}{1-e^{i\theta}}\right)-\frac{1}{n_k}\left(\frac{1-e^{in_k\theta}}{1-e^{i\theta}}\right)\right|\\
&=\frac{1}{|1-e^{i\theta}|}\cdot \left|\frac{1-e^{im\theta}}{\theta m}-\frac{1-e^{in_k\theta}}{\theta n_k}\right|\\
&\leq 4 \left|\frac{1-e^{im\theta}}{\theta m}-\frac{1-e^{in_k\theta}}{\theta n_k}\right|\\
&\leq \frac{8}{\theta m}+\frac{8}{\theta n_k}\\
&\leq \frac{16}{\theta n_k}
\end{align*}
since $n_k\leq m$, and this estimate allows us to control $J_2(\theta )$ as the way we controlled $I_2(\theta )$  in the proof of Theorem~\ref{varineq}. $J_1(\theta )$ can also be controlled as the way we controlled $I_1(\theta )$ in the  in the proof of Theorem~\ref{varineq}. It is clear that 
\begin{align*}
J_1(\theta )&=\sum_{k:n_k\theta <1} \left(\int_{-\pi}^{\pi}\sup_{\substack{n_k\leq m< n_{k+1}\\m\in \mathbb{N}}}|a_m(\theta )-a_{n_k}(\theta )|^2\, d\mu (\theta )\right)^{1/2}\\
&\leq 4\sum_{k:n_k\theta <1} \left(\int_{-\pi}^{\pi}\sup_{\substack{n_k\leq m< n_{k+1}\\m\in \mathbb{N}}}\left|\frac{1-e^{im\theta}}{\theta m}-\frac{1-e^{in_k\theta}}{\theta n_k}\right|^2\, d\mu (\theta )\right)^{1/2}.
\end{align*}
\indent
Let $F:\mathbb{R}\to\mathbb{C}$ be the same function as in the proof of Theorem~\ref{varineq}. By the Mean Value Theorem there exits a constant $\xi\in (n_k, m)$ such that
$$|F(m )-F(n_k)|= |F^\prime (\xi )| |m-n_k|.$$
\indent
Since $\xi \theta \leq m\theta <n_{k+1}\theta <1$, we have
\begin{align*}
|F(m)-F(n_k)|&= |F^\prime (\xi )| |m-n_k|\\
&\leq C\theta (m-n_k)\\
&\leq C\theta m, \;\;\;\; ({\text{since}}\; \;0 <n_k\leq m) \\
&< C\theta n_{k+1}\cdot \frac{n_{k_0-1}}{n_{k_0-1}}\\
&=C\frac{n_{k+1}}{n_{k_0-1}}\cdot n_{k_0-1}\theta \\
&<C\frac{n_{k+1}}{n_{k_0-1}}, \;\;\;\; ({\text{since}}\; \;n_{k_0-1}\theta <1),
\end{align*}
and rest follows as the way we controlled $I_1(\theta )$ in Theorem~\ref{varineq}.
\end{proof}
\begin{cor}\label{ercorosc}Let  $(X,\mathscr{B} ,\mu ,\tau )$ be an ergodic, measure preserving dynamical system and define the ergodic averages
$$\mathcal{A}_nf(x)=\frac{1}{n}\sum_{j=1}^nf\circ\tau^j(x).$$
Let $(n_k)$  be a lacunary sequence, and let $\mathbb{N}$ be the set of natural numbers. Then there exists a positive constant $C$ such that
$$\sum_{k=1}^\infty\sup_{\substack{n_k\leq m< n_{k+1}\\m\in \mathbb{N}}}\|\mathcal{A}_mf-\mathcal{A}_{n_k}f\|_2\leq C\|f\|_2$$
for all $f\in L^2(X)$.
\end{cor}
\begin{proof}It is clear  that $Tf(x)=f\circ\tau (x)$ is a unitary operator,  and thus it is a contraction. Since $L^2(X)$ is a Hilbert space, the result follows from Theorem~\ref{osineq}.
\end{proof}
\begin{rem} Since a unitary operator $T$ on a Hilbert space is a contraction, our results remain also true for a unitary operator.  On the other hand,  one can give a direct proof  when $T$ is a unitary operator in Theorem~\ref{varineq} and Theorem~\ref{osineq} by using the spectral theorem for a unitary operator on a complex Hilbert space (see F.~Riesz and B.~Sz-Nagy~\cite{frbsn}, p. 289). After proving these theorems for a unitary operator $T$ we can also prove them for a contraction as follows:\\
\indent
By the dilation theorem (see B.~Sz-Nagy and C.~Foias~\cite{bsncf}, p. 45), there exists a Hilbert space $\mathcal{L}$ containing $\mathcal{H}$ as a closed subspace, an orthogonal projection $P:\mathcal{L}\to \mathcal{H}$, and a unitary operator $U:\mathcal{L}\to \mathcal{L}$ with $PU^jf=T^jf$ for all $j\geq 0$ and $f\in \mathcal{H}$. But then for $f\in \mathcal{H}$, and for any $N\in \mathbb{Z}^+$ we obtain
\begin{align*}
\sum_{k=1}^N\|A_{n_{k+1}}(T)f-A_{n_k}(T)f\|_{\mathcal{H}}&=\sum_{k=1}^N\|P(A_{n_{k+1}}(U)f-A_{n_k}(U)f)\|_{\mathcal{H}}\\
&=\|P\|\sum_{k=1}^N\|A_{n_{k+1}}(U)f-A_{n_k}(U)f\|_{\mathcal{H}}\\
&\leq C
\end{align*}
since $\|P\|\leq 1$. This proves Theorem~\ref{varineq}, and Theorem~\ref{osineq} follows from the following inequality similarly:
\begin{align*}
\sum_{k=1}^N \sup_{\substack{n_k\leq m< n_{k+1}\\m\in \mathbb{N}}} \|A_m(T)f-A_{n_k}(T)f\|_{\mathcal{H}}&=\sum_{k=1}^N\sup_{\substack{n_k\leq m< n_{k+1}\\m\in \mathbb{N}}}\|P(A_m(U)f-A_{n_k}(U)f)\|_{\mathcal{H}}\\
&=\sum_{k=1}^N\sup_{\substack{n_k\leq m< n_{k+1}\\m\in \mathbb{N}}}\|P\|\| A_m(U)f-A_{n_k}(U)f  \|_{\mathcal{H}}\\
&=\|P\|\sum_{k=1}^N\sup_{\substack{n_k\leq m< n_{k+1}\\m\in \mathbb{N}}}\|A_m(U)f-A_{n_k}(U)f\|_{\mathcal{H}}\\
&\leq C.
\end{align*}
\end{rem}

\end{document}